\documentclass[12pt]{amsart}
\usepackage{amsmath}
\usepackage{amssymb} 
\usepackage{amscd}
\usepackage[mathscr]{eucal} 
\usepackage{latexsym}

\newtheorem{thm}{Theorem}[section]

\theoremstyle{definition}

\newtheorem{defn}[thm]{Definition}

\newtheorem{ex}[thm]{Example}

\theoremstyle{remark}
\newtheorem{notation}{Notation}

\numberwithin{equation}{section}

\def\F{{\mathbb F}}

\def\Ext{\text{\rm Ext}}

\def\hom{\text{\rm hom}}

\def\F{\text{\rm F}}
\def\GMHS{\text{\rm GMHS}}
\def\Gr{\text{\rm Gr}}

\def\Im{\text{\rm Im}\,}
\def\Ker{\text{\rm Ker}\,}

\def\MHS{\text{\rm MHS}}
\def\CH{\text{\rm CH}}

\def\Spec{\text{\rm Spec}\,}

\def\lim{\text{\rm lim}}

\begin{document} 

\title[Generalization of the theory of mixed Hodge structures]
{Generalization of the theory of mixed Hodge structures and its application}

\author[Kazuma Morita]{Kazuma Morita}
\address{Department of Mathematics, Hokkaido University, Sapporo 060-0810, Japan}
\email{morita@math.sci.hokudai.ac.jp}

\subjclass{ 
14C15, 14C30
 }
\keywords{ 
 Hodge theory, Cycles, Chow groups}
\date{\today}

\maketitle
\begin{abstract}In this paper, we shall generalize the theory of mixed Hodge structures due to Deligne and obtain a subcategory GMHS in the category of mixed Hodge structures such that we have $\Ext_{\text{GMHS}}^{2}(\mathbb{Q},-)\not=0$ in general. 
\end{abstract}
\section{Introduction}
For a smooth projective scheme $X$ over $\mathbb{C}$, there exists a conjectural filtration $\F_{M}$ (called BBM filtration) on the Chow group $\CH^{r}(X, \mathbb{Q})=\CH^{r}(X)\otimes\mathbb{Q}$ such that we have
$\Gr_{\F_{M}}^{m}\CH^{r}(X, \mathbb{Q})=\Ext_{M}^{m}(\mathbb{Q},H^{2r-m}(X)(r))$. Here, $M$ is the conjectural category of mixed motives over $\Spec(\mathbb{C})$. On the other hand, by the realization functor $M\rightarrow \MHS$ from the category of mixed motives to  that of mixed Hodge structures, we should have $$\Gr_{\F_{M}}^{m}\CH^{r}(X, \mathbb{Q})\rightarrow \Ext_{\MHS}^{m}(\mathbb{Q},H^{2r-m}(X(\mathbb{C}),\mathbb{Q}(r))).$$
It is well-known, however, that the higher extension group $\Ext_{\MHS}^{m}(\mathbb{Q},-)$ for $2\leq m$ always vanishes and one cannot obtain any information about the graded piece $\Gr_{\F_{M}}^{m}\CH^{r}(X, \mathbb{Q})$ for $2\leq m$ by using the extension of mixed Hodge structures.

In this paper, we shall generalize the theory of mixed Hodge structures due to Deligne and obtain a subcategory GMHS in the category of mixed Hodge structures such that we have $\Ext_{\text{GMHS}}^{2}(\mathbb{Q},-)\not=0$ in general. Note that M.Asakura constructs another category (called the category of arithmetic Hodge structures) and shows that the higher extension group does not vanish in this category [A]. One will see that the category GMHS is an abelian category and that there is a forgetful functor $\mathscr{F}$ to the category MHS. 

{\bf Acknowledgments} 
The author would like to thank his advisor Professor
Kazuya Kato for continuous advice, encouragements and patience. 
He is also grateful to Professor Masanori Asakura for useful discussions.
A part of this work was done while he was staying at Universit\'e Paris-Sud 11 and
he thanks this institute for the hospitality.
His staying at Universit\'e Paris-Sud 11 was partially 
supported by JSPS Core-to-Core Program
``New Developments of Arithmetic Geometry, Motives, Galois Theory, and Their Practical Applications''
and he thanks Professor Makoto Matsumoto for encouraging this visiting.   
This research was partially supported by JSPS Research Fellowships for Young Scientists.
\section{Mixed Hodge structures}
\subsection{Review of the classical theory}
For a compact K\"ahler manifold $X$, Hodge shows that there exists a decomposition
$$H^{k}(X,\mathbb{C})=\oplus_{p+q=k}H^{p,q}(X)$$
where $H^{p,q}(X)$ is a complex subspace and satisfies the Hodge symmetry $H^{p,q}(X)=\overline{H^{q,p}(X)}$ ($-$ denotes the complex conjugation). This leads to the following definition.
\begin{defn}
An integral Hodge structure of pure weight $k$ is a free abelian group $H_{\mathbb{Z}}$ of finite type equipped with a decomposition 
$$H_{\mathbb{C}}:=H_{\mathbb{Z}}\otimes \mathbb{C}=\oplus_{p+q=k}H^{p,q}$$
where  $H^{p,q}$ is a complex subspace and satisfies the symmetry $H^{p,q}=\overline{H^{q,p}}$.
\end{defn}
Let $H_{\mathbb{Z}}$ be an integral Hodge structure of pure weight $k$ and define a decreasing filtration (called Hodge Filtration) $F^{\cdot}H_{\mathbb{C}}$ by $F^{p}H_{\mathbb{C}}=\oplus_{p\leq r}H^{r,k-r}$. This filtration satisfies $H_{\mathbb{C}}=\F^{p}H_{\mathbb{C}}\oplus \overline{F^{k-p+1}H_{\mathbb{C}}}$ and determines the Hodge decomposition by the formula $H^{p,q}=F^{p}H_{\mathbb{C}}\cap \overline{ F^{q}H_{\mathbb{C}}}$. Let HS be the category of Hodge structures of pure weights: its object is given by Hodge structure $(H_{\mathbb{Z}}, F^{\cdot}H_{\mathbb{C}})$ of pure weight and its morphism is given by a morphism $f:H_{\mathbb{Z}}\rightarrow H_{\mathbb{Z}}'$ which is compatible with the Hodge filtration $F^{\cdot}$.

For a general separated scheme $X$ of finite type over $\mathbb{C}$, the cohomology group $H^{k}(X,\mathbb{Z})$ does not carry the Hodge structure of pure weight in general. Then, Deligne shows that there exists an increasing filtration $W_{\cdot}$ (called weight filtration) on $H^{k}(X,\mathbb{Z})$ such that the Hodge filtration induced on  $\Gr_{r}^{W}H^{k}(X,\mathbb{C})$ defines an integral Hodge structure of pure weight $k+r$ on 
$\Gr_{r}^{W}H^{k}(X,\mathbb{Z})$ ([D1], [D2]). This result leads to the following definition.
\begin{defn}
A mixed Hodge structure of weight $k$ is a free abelian group $H_{\mathbb{Z}}$ of finite type equipped with an increasing filtration (called weight filtration) $W_{\cdot}$ on $H_{\mathbb{Z}}$ and a decreasing filtration (called Hodge filtration) $F^{\cdot}$ on $H_{\mathbb{C}}$ such that  the filtration induced by $F^{\cdot}$ on $\Gr_{r}^{W}H_{\mathbb{C}}$ defines a Hodge structure of pure weight $k+r$ on $\Gr_{r}^{W}H_{\mathbb{Z}}$.
\end{defn}
Let MHS be the category of mixed Hodge structures: its object is given by a mixed Hodge structure $(H_{\mathbb{Z}}, F^{\cdot}H_{\mathbb{C}}, W_{\cdot}H_{\mathbb{Z}})$ and its morphism is given by a morphism $f:H_{\mathbb{Z}}\rightarrow H_{\mathbb{Z}}'$ which is compatible with filtrations $F^{\cdot}$ and $W_{\cdot}$.
It is known that the category MHS is an abelian category ([D1], (2.3.5)).
\subsection{Generalized mixed Hodge structures}
In this subsection, we shall generalize the theory of mixed Hodge structures due to Deligne. Let $U$ be a smooth and separated scheme of finite type over $\mathbb{C}$ and $X$ be a smooth compactification of $U$.  By using subschemes on $U$ and $D=X\verb+\+U$, we shall introduce two structures on cohomology groups: $z$-structures and $w$-structures. These lead to the generalization of the theory of mixed Hodge structures.  
\subsubsection{$z$-structures on cohomology groups}
For a subscheme $V$ on $U$, let  $\vec{r}=(r_{1},\cdots,r_{l})$ denote a basis of 
$\Im(H_{V}^{n}(U,\mathbb{Q})\rightarrow H^{n}(U,\mathbb{Q}))$ over $\mathbb{Q}$. Choose elements   
$\vec{s}=(s_{1},\cdots,s_{m})$ of $H^{n}(U,\mathbb{Q})$ such that  $\verb+{+\vec{r},\vec{s}\verb+}+$ forms 
a basis of $H^{n}(U,\mathbb{Q})$ over $\mathbb{Q}$. 
Then, define an involution $z_{V,\vec{r},\vec{s}}$ on $H^{n}(U,\mathbb{Q})$ by the formula 
$$z_{V,\vec{r},\vec{s}}(\sum_{i=1}^{l}a_{i}r_{i}+\sum_{j=1}^{m} b_{j}s_{j})=-\sum_{i=1}^{l} a_{i}r_{i}+\sum_{j=1}^{m} b_{j}s_{j}\quad (a_{i},b_{j}\in\mathbb{Q}).$$
\subsubsection{$w$-structures on cohomology groups}
Since we assume that $U$ is a smooth and separated scheme of finite type over $\mathbb{C}$, it is a Zariski open set in a complete scheme $X$ ([N]). Furthermore, we assume that $X$ is smooth projective and that the complement $D=X\verb+\+U$ is a globally normal crossing divisor, that is, we have $D=\cup_{i\in I} D_{i}$ where each $D_{i}\subset X$ is a smooth hypersurface and the intersection of hypersurfaces is transverse ([H]).  
\begin{notation}
For a subset $K\subset I$, put $D_{K}=\cap_{i\in K}D_{i}$ and let $D^{(k)}$ denote the disjoint union of $D_{K}$ where $K$ runs through subsets of $I$ of cardinal $k$. Set $D^{(0)}=X$. 
\end{notation}  
For the weight spectral sequence ${}_W E$ associated to the weight filtration $W_{\cdot}$, we have 
${}_W E_{1}^{p,q}\simeq H^{2p+q}(D^{(-p)},\mathbb{C})$ and its differential $d_{1}$ is given by
\begin{equation}
\begin{CD}
H^{2p+q}(D^{(-p)},\mathbb{C})    @>\text{$d_{1}$}>>   H^{2p+q+2}(D^{(-p-1)},\mathbb{C}) \\
@|    @|  \\
\bigoplus _{|K|=-p}H^{2p+q}(D_{K},\mathbb{C}) @>\text{$d_{1}$}>> \bigoplus _{|L|=-p-1}H^{2p+q+2}(D_{L},\mathbb{C})
\end{CD}
\end{equation}
where $d_{1}$ has the component $d_{1K}^{\mspace{10mu}L}$ equal to zero for $L\not\subset K$ and equal to $(-1)^{q+s}j_{K*}^{L}$ for $K=\verb+{+i_{1}<\cdots<i_{p}\verb+}+$ and $L=K\verb+\+\verb+{+i_{s}\verb+}+$ where $j_{K*}^{L}$ denotes the Gysin map corresponding to the inclusion $j_{K}^{L}:D_{K}\hookrightarrow D_{L}$. Due to the result of Deligne, this spectral sequence degenerates at $E_{2}$ and we obtain ${}_W E_{2}^{p,q}=\Gr_{-p}^{W}H^{k}(U,\mathbb{C})$.  
For each subscheme $V'$ on $D$, we shall define a $\mathbb{C}$-linear involution $w_{V'}$ on $\Gr_{-p}^{W}H^{k}(U,\mathbb{C})$. If we have $V'\not\subset D^{(-p)}$, put $w_{V'}(c)=c$. 
Now, assume that $V'$ is a subscheme on $D^{(-p)}$ and then there is a natural morphism
$$\psi_{V'}:H^{2p+q}_{V'}(D^{(-p)},\mathbb{C})\rightarrow  H^{2p+q}(D^{(-p)},\mathbb{C})\simeq {}_W E_{1}^{p,q} \twoheadrightarrow{}_W E_{1}^{p,q}/\Im({}_W E_{1}^{p-1,q}).$$
Let  $\vec{r}=(r_{1},\cdots,r_{l})$ denote a basis of $\Im(\psi_{V'})\cap {}_W E_{2}^{p,q}$ over $\mathbb{C}$. Choose elements   
$\vec{s}=(s_{1},\cdots,s_{m})$ of ${}_W E_{2}^{p,q}$ such that  $\verb+{+\vec{r},\vec{s}\verb+}+$ forms 
a basis of ${}_W E_{2}^{p,q}$ over $\mathbb{C}$. 
Then, define an involution $w_{V',\vec{r},\vec{s}}$ on $\Gr_{-p}^{W}H^{k}(U,\mathbb{C})$  by the formula 
$$w_{V',\vec{r},\vec{s}}(\sum_{i=1}^{l}a_{i}r_{i}+\sum_{j=1}^{m} b_{j}s_{j})=-\sum_{i=1}^{l} a_{i}r_{i}+\sum_{j=1}^{m} b_{j}s_{j}\quad (a_{i},b_{j}\in\mathbb{C}).$$
\subsubsection{Category of generalized mixed Hodge structures}
The results of preceding subsections lead to the following definition.
\begin{defn}Let $U$ be a smooth and separated scheme of finite type over $\mathbb{C}$ and $X$ be a smooth compactification of $U$ such that $D=X\verb+\+U$ is a globally normal crossing divisor. 
A generalized mixed Hodge structure consists of $(H_{\mathbb{Z}},F^{\cdot},W_{\cdot},\verb+{+z_{V}\verb+}+_{V\subset U},\verb+{+w_{V'}\verb+}+_{V'\subset D})$ where 
\begin{itemize}
\item the triple $(H_{\mathbb{Z}},F^{.},W_{\cdot})$ is a mixed Hodge structure, 
\item $z_{V}$ denotes a $\mathbb{C}$-linear involution on $H_{\mathbb{C}}$ for each subscheme $V$ on $U$,
\item $w_{V'}$ denotes a $\mathbb{C}$-linear isomorphism of $H_{\mathbb{C}}$ such that the induced action on  $\Gr_{m}^{W}H_{\mathbb{C}}$ is an involution for each subscheme $V'$ on $D$.
\end{itemize}
\end{defn}
Let GMHS denote the category of generalized mixed Hodge structures: its object is given by a generalized mixed Hodge   structure and its morphism between                          $\bigl\{(H_{\mathbb{Z}}^{i},F^{\cdot},W_{\cdot},\verb+{+z_{V_{i}}\verb+}+_{V_{i}\subset  U_{i}},\verb+{+w_{V'_{i}}\verb+}+_{V'_{i}\subset D_{i}})\bigr\}_{i=1,2}$  is given by the pair of a morphism of mixed Hodge structures $f:H_{\mathbb{Z}}^{1}\rightarrow H_{\mathbb{Z}}^{2}$ and a morphism of schemes $g:X_{2}\rightarrow X_{1}$ where $X_{i}$ denotes a smooth compactification of $U_{i}$ such that we have $D_{i}=X_{i}\verb+\+U_{i}$. Furthermore,  assume that this pair of morphisms satisfies the compatible condition $f\circ x=y\circ f$ where 
$$
     \left\{ 
     \begin{aligned}
x=&z_{V_{1}}, \quad \mspace{10mu}y=z_{V_{2}} \qquad \mspace{6mu}\text{if $g(U_{2}-V_{2})\mspace{3mu}\subset U_{1}-V_{1}$ $\mspace{7mu}$and $g:V_{2}\simeq V_{1}$},\\   
x=&z_{V_{1}}, \quad \mspace{10mu}y=w_{V_{2}'} \qquad \mspace{1mu}\text{if $g(D_{2}-V_{2}')\subset U_{1}-V_{1} $ $\mspace{6mu}$and $g:V_{2}'\simeq V_{1}$}, \\
x=&w_{V_{1}'}, \quad \mspace{6mu}y=z_{V_{2}} \qquad \mspace{5mu}\text{if $g(U_{2}-V_{2})\mspace{3mu}\subset D_{1}-V_{1}'$ $\mspace{4mu}$and $g:V_{2}\simeq V_{1}'$},\\
x=&w_{V_{1}'}, \quad \ y=w_{V_{2}'} \qquad \text{if $g(D_{2}-V_{2}')\subset D_{1}-V_{1}' \mspace{3mu}$ and $g:V_{2}'\simeq V_{1}'$}.
     \end{aligned}
     \right.$$
One can verify that the category GMHS is an abelian category and that there is a forgetful functor $\mathscr{F}$ to the category MHS.
\section{Extension  groups $\Ext^{m}$}
For a smooth projective scheme $X$ over $\mathbb{C}$, the conjectural filtration $\F_{M}$ on the Chow group $\CH^{r}(X, \mathbb{Q})$ should satisfy
$\Gr_{\F_{M}}^{m}\CH^{r}(X, \mathbb{Q})=\Ext_{M}^{m}(\mathbb{Q},H^{2r-m}(X)(r))$. Here, $M$ is the conjectural category of mixed motives over $\Spec(\mathbb{C})$. On the other hand, by the realization functor $M\rightarrow \MHS$, we should have $$\Gr_{\F_{M}}^{m}\CH^{r}(X, \mathbb{Q})\rightarrow \Ext_{\MHS}^{m}(\mathbb{Q},H^{2r-m}(X(\mathbb{C}),\mathbb{Q}(r))).$$
From the right exactness of $\Ext^{1}_{\text{MHS}}(\mathbb{Q},-)$, however, it follows that the higher extension group $\Ext_{\MHS}^{m}(\mathbb{Q},-)$ for $2\leq m$ always vanishes [C]. In this section, we shall introduce the one dimensional vector space $\mathbb{Q}_{M}$ over $\mathbb{Q}$ equipped with generalized mixed Hodge structures and shall construct an example which shows that $\Ext_{\text{GMHS}}^{2}(\mathbb{Q}_{M},-)$ does not vanish in general.
\subsection{Definition of $\mathbb{Q}_{M}$}
Let $U$ be a smooth and separated scheme of finite type over $\mathbb{C}$ and $X$ be a smooth compactification of $U$   such that $D=X\verb+\+ U$ is a globally normal crossing divisor. Define $S(U,D)$  to be the set of subschemes $(\verb+{+V\verb+}+_{V\subset U},\verb+{+V'\verb+}+_{V'\subset D})$ where $V$ (resp. $V'$) runs through any subscheme on $U$ (resp. $D$).
\begin{defn}
With notations as above, for a subset $M(U,D)$ of $S(U,D)$, let $\mathbb{Q}_{M(U,D)}$ denote the one  dimensional vector space over $\mathbb{Q}$ equipped with the generalized mixed Hodge structure $(\mathbb{Q},F^{\cdot},W_{\cdot},$ $\verb+{+z_{V}\verb+}+_{V\subset U},\verb+{+w_{V'}\verb+}+_{V'\subset D})$ where the mixed Hodge structure is trivial and involutions $(\verb+{+z_{V}\verb+}+_{V\subset U},\verb+{+w_{V'}\verb+}+_{V'\subset D})$ act on $\mathbb{Q}_{M(U,D)}$ by 
$$
     \left\{
     \begin{aligned}
          z_{V}(a)&=-a \ \text{if $V\in M(U,D)$},\quad \mspace{11mu}z_{V}(a)=a  \ \text{if $V\not\in M(U,D)$} \\
          w_{V'}(a)&=-a \ \text{if $V'\in M(U,D)$},\quad  w_{V'}(a)=a  \ \text{if $V'\not\in M(U,D)$}
     \end{aligned}
     \right.$$
\end{defn}
\begin{ex}
Let $X$ be a smooth projective scheme over $\mathbb{C}$ and let $cl:\CH^{j}(X,\mathbb{Q})$ $ \rightarrow H^{2j}(X(\mathbb{C}),\mathbb{C})$ denote the cycle map. Then, the classical Hodge conjecture states that this cycle map has the image 
$$H^{j,j}(X) \cap H^{2j}(X(\mathbb{C}),\mathbb{Q})=\Ext^{0}_{\text{HS}}(\mathbb{Q}, H^{2j}(X(\mathbb{C}),\mathbb{Q})).$$ 
Assume that the classical Hodge conjecture holds. Then, we can write $\Ext_{\text{HS}}^{0}(\mathbb{Q},$ $ H^{2j}(X(\mathbb{C}),\mathbb{Q}))$ in terms of generalized mixed Hodge structures 
$$\Ext^{0}_{\text{HS}}(\mathbb{Q}, H^{2j}(X(\mathbb{C}),\mathbb{Q}))=\bigoplus_{M\subset S(X,\phi) }\Ext^{0}_{\text{GMHS}}(\mathbb{Q}_{M}, H^{2j}(X(\mathbb{C}),\mathbb{Q})).$$ 
\begin{proof}
It suffices to show that  we have LHS $\subset$ RHS.  
Note that the cohomology group $H^{2j}(X(\mathbb{C}),\mathbb{Q})$ is equipped with the involution $z_{V}$ for each subscheme $V$ on $X$ through $H^{2j}_{V}(X,\mathbb{Q})\rightarrow H^{2j}(X,\mathbb{Q})$. 
Take an element $f_{a} \ (:1\mapsto a)$ of $\Ext^{0}_{\text{HS}}(\mathbb{Q}, H^{2j}(X(\mathbb{C}),\mathbb{Q}))$. By the assumption, there is an element $\widetilde{a}$ of $\CH(X,\mathbb{Q})$ such that we have $f_{a}=cl(\widetilde{a})$. 
We can write this cycle $\widetilde{a}$ as  $\sum_{k=1}^{m}n_{k}\widetilde{a_{k}}$ ($n_{k}\in\mathbb{Z}$) where $\verb+{+\widetilde{a_{k}}\verb+}+_{k=1}^{m}$ denote subschemes on $X$.  
Then, the element $f_{a}=\sum_{k=1}^{m}n_{k}cl(\widetilde{a}_{k})$ is contained in $\bigoplus_{k=1}^{m}(\bigoplus_{M_{k}} \Ext^{0}_{\text{GMHS}}(\mathbb{Q}_{M_{k}}, H^{2j}(X(\mathbb{C}),\mathbb{Q})))$ where $M_{k}$ runs through any set containing  $\bigr(\verb+{+\widetilde{a}_{k}\verb+}+, \verb+{+\phi\verb+}+\bigl)$ and thus we obtain LHS $\subset$ RHS.
\end{proof}
\end{ex}
\begin{ex}
With notations as in the previous example, define $\CH^{j}(X,\mathbb{Q})_{\hom}=\verb+{+\alpha\in \CH^{j}(X,\mathbb{Q})\verb+|+ cl(\alpha)=0\verb+}+$. Then, we have the Abel-Jacobi map
$$cl':\CH^{j}(X,\mathbb{Q})_{\hom}\rightarrow \frac{H^{2j-1}(X(\mathbb{C}),\mathbb{C})}{F^{j}H^{2j-1}(X(\mathbb{C}),\mathbb{C})\oplus H^{2j-1}(X(\mathbb{C}),\mathbb{Q})}.$$
One can see that the target of this map is isomorphic to $\Ext_{\text{MHS}}^{1}(\mathbb{Q},H^{2j-1}(X,\mathbb{Q}))$. First, we shall review the construction of the extension class given by the Abel-Jacobi map. Let $v$ denote an element of $\CH^{j}(X,\mathbb{Q})_{\hom}$ and $V$ be the support of $v$. Put $U=X\verb+\+V$. Then, there exists a commutative diagram 
$$
\begin{CD}
0 @>>> H^{2j-1}(X,\mathbb{Q})@>>> H^{2j-1}(U,\mathbb{Q})@>>> H^{2j}_{V}(X,\mathbb{Q})@>>> H^{2j}(X,\mathbb{Q})\\
@. @|          @A\cup AA        @A h_{v} AA                 @. \\
0 @>>>H^{2j-1}(X,\mathbb{Q})@>>> E @>>> \mathbb{Q} @>>> 0
\end{CD}
$$
where $h_{v}$ denotes the map $\mathbb{Q}\rightarrow H^{2j}_{V}(X,\mathbb{Q}):1\mapsto v$ and the bottom exact sequence is obtained by pull-back via $h_{v}$. We can verify that the extension class $E$ of this exact sequence is the image of $v$ under the Abel-Jacobi map ([J1], 9.4).    
Now, let us see the bottom exact sequence of the diagram above
in terms of generalized mixed Hodge structures. 
We can write the cycle $v$ as $\sum_{k=1}^{m}n_{k}v_{k}$ ($n_{k}\in\mathbb{Z}$) where
$\verb+{+v_{k}\verb+}+_{k=1}^{m}$ denote subschemes on $X$. We will denote $\mathbb{Q}$ in the diagram above by $\mathbb{Q}_{v}$ and fix a basis $1_{v}$ of $\mathbb{Q}_{v}$ over $\mathbb{Q}$.  Since the involution $v\mapsto -v$ on $H_{V}^{2j}(X,\mathbb{Q})$ should correspond to the involution $f_{v}:1_{v}\mapsto -1_{v}$, it is natural to think that $\mathbb{Q}_{v}$ is contained in $\mathbb{Q}_{N}=\bigoplus_{k=1}^{m}( \bigoplus_{M_{k}}\mathbb{Q}_{M_{k}})$ where $M_{k}$ runs through any set of $S(X,\phi)$ containing $(\verb+{+v_{k}\verb+}+,\verb+{+\phi\verb+}+)$.  Since the extension class $E$ is the image of $v$  under the Abel-Jacobi map, it is compatible with the action of $f_{v}$ induced by $z$-structures on $\mathbb{Q}_{N}$. Furthermore, $E$ is clearly compatible with $z$-structures on $H^{2j-1}(X,\mathbb{Q})$. Thus, we can regard the extension class $E$ as an element of    
$\bigoplus_{M\subset S(X,\phi)}\Ext_{\text{GMHS}}^{1}(\mathbb{Q}_{M},H^{2j-1}(X,\mathbb{Q}))$. 
\end{ex}
\subsection{Non-vanishing of $\Ext^{2}_{\GMHS}(\mathbb{Q}_{M},-)$}
In this section, we shall see that the higher extension group $\Ext^{2}_{\text{GMHS}}(\mathbb{Q}_{M},-)$ does not vanish in general. First, we shall recall the Yoneda extension class. Let $A$ denote an abelian category. For objects $M$ and $N$ of $A$, an element $C(E)$ of $\Ext_{A}^{n}(M,N)$ (called the Yoneda extension class) is given by an exact sequence 
$$E:0\rightarrow N\rightarrow R_{n}\rightarrow R_{n-1}\rightarrow \cdots \rightarrow R_{1}\rightarrow M \rightarrow 0.$$   
Let $E':0\rightarrow N\rightarrow R_{n}'\rightarrow R_{n-1}'\rightarrow \cdots \rightarrow R_{1}'\rightarrow M \rightarrow 0$ be another extension. Then, we have $C(E)=C(E')$ if and only if there exists an extension $E'':0\rightarrow N\rightarrow R_{n}''\rightarrow R_{n-1}''\rightarrow \cdots \rightarrow R_{1}''\rightarrow M \rightarrow 0$ such that we have the following commutative diagram
$$
\begin{CD}
0 @>>> N   @>>> R_{n}   @>>> R_{n-1}  @>>> \cdots @>>> R_{1}  @>>>    M  @>>> 0      \\
@.     @|       @AAA          @AAA            @.         @AAA         @|      @.     \\
0 @>>> N   @>>> R_{n}''  @>>> R_{n-1}'' @>>> \cdots @>>> R_{1}'' @>>> M  @>>> 0      \\
@.     @|       @VVV          @VVV            @.        @VVV          @|      @.     \\
0 @>>> N   @>>> R_{n}'  @>>> R_{n-1}' @>>> \cdots @>>> R_{1}' @>>>    M  @>>> 0.      \\
\end{CD}
$$
\begin{ex}
We shall construct an example which shows that $\Ext^{2}_{\GMHS}(\mathbb{Q}_{M},-)$ does not vanish.
Let us consider the following exact sequence in the category of generalized mixed Hodge structures 
$$E: 0\rightarrow S \xrightarrow{i} T \xrightarrow{j} V \xrightarrow{k} \mathbb{Q}_{M}\rightarrow 0.$$
Here, 
\begin{itemize}
\item  $S$ is a $2$-dimensional vector space over $\mathbb{Q}$ equipped with the Hodge structure of pure weight $-1$. For a smooth projective curve $X$ over $\mathbb{C}$, assume that $S$ is endowed with the trivial $z$-structure on $X$ and the trivial $w$-structure on $\phi \ (=X\verb+\+X)$.\\
\item  $T$ is a $3$-dimensional vector space over $\mathbb{Q}$ equipped with the mixed Hodge structure such that 
$\Gr_{0}^{W}(T)$ has the Hodge structure of pure weight $-1$ and $\Gr_{1}^{W}(T)$ has the Hodge structure of pure weight $0$. For two points $\verb+{+ D_{i}\verb+}+_{i=1,2}$ on $X$, assume that $T$ is equipped with the trivial $z$-structure on $U=X\verb+\+\verb+{+D_{i}\verb+}+_{i=1,2}$ and with the $w$-structure on $\verb+{+D_{i}\verb+}+_{i=1,2}$. Here, the action of $\verb+{+w_{D_{i}}\verb+}+_{i=1,2}$ on a basis $\verb+{+ e,s_{1},s_{2}\verb+}+$ of $T$ over $\mathbb{Q}$ is given by
$$w_{D_{i}}
\begin{pmatrix}
 e\\
 s_{1}\\
 s_{2}\\
\end{pmatrix}
=
\begin{pmatrix}
1&-1&0\\
0&1&0\\
0&0&1
\end{pmatrix}
\begin{pmatrix}
 e\\
 s_{1}\\
 s_{2}\\
\end{pmatrix}
$$ 
where $\verb+{+ s_{1},s_{2}\verb+}+$ denotes a basis of $i(S)$ over $\mathbb{Q}$.
Note that these actions induce (trivial) involutions on  $\Gr_{i}^{W}(T)$ ($i=0,1$).\\
\item $V$ is a $2$-dimensional vector space over $\mathbb{Q}$ equipped with the Hodge structure of pure weight $0$. 
Assume that $V$ is equipped with the trivial $w$-structure on $\phi$ and with the $z$-structure on $\verb+{+D_{i}\verb+}+_{i=1,2}$ such that the action of $\verb+{+z_{D_{i}}\verb+}+_{i=1,2}$ on a basis $\verb+{+\alpha,\beta\verb+}+$ of $V$ over $\mathbb{Q}$ is given by
$$
z_{D_{i}}
\begin{pmatrix}
\alpha \\
\beta
\end{pmatrix}
=
\begin{pmatrix}
1 & 0 \\
c(D_{i}) & -1
\end{pmatrix}
\begin{pmatrix}
\alpha \\
\beta
\end{pmatrix}
$$ 
where $\verb+{+c(D_{i})\verb+}+_{i=1,2}$  satisfy $c(D_{1})\not=c(D_{2})$. For  a subscheme $D'$ on $X$ other than $\verb+{+D_{i}\verb+}+_{i=1,2}$, assume that the action of $z_{D'}$ on $V$  is trivial. Then, we can consider that $V$ is also equipped with the $z$-structure on $X$ and the trivial $w$-structure on $\phi$.\\
\item $M=(\verb+{+D_{i}\verb+}+_{i=1,2},\verb+{+\phi\verb+}+)\subset S(X,\phi)$, that is, $\mathbb{Q}_{M}$ is endowed with the  non-trivial action of $z_{D_{i}}$ and  the trivial $w$-structure on $\phi$.
\end{itemize}
On the other hand, one can verify that the exact sequence $E': 0\rightarrow S \rightarrow S \rightarrow  \mathbb{Q}_{M} \rightarrow \mathbb{Q}_{M}\rightarrow 0$ in GMHS gives a trivial Yoneda extension class. 
Thus, it suffices to show that we have $C(E)\not=C(E')$ in $\Ext^{2}_{\GMHS}(\mathbb{Q}_{M},S)$. Assume that there exists an exact sequence $0\rightarrow S \rightarrow T' \rightarrow V' \rightarrow \mathbb{Q}_{M}\rightarrow 0$ in GMHS such that we have the following commutative diagram
\begin{align}
\begin{CD}
0 @>>> S   @> i >> T    @> j >> V   @> k >> \mathbb{Q}_{M} @>>> 0        \\
@.     @|       @AA f A    @AA g A            @|             @.       \\
0 @>>> S   @>>> T'  @> h> > V' @> h' >> \mathbb{Q}_{M} @>>> 0        \\
@.     @|       @VV f' V    @VV g' V            @|             @.       \\
0 @>>> S   @>>> S   @>>>   \mathbb{Q}_{M} @>>> \mathbb{Q}_{M} @>>> 0.       \\
\end{CD}
\end{align}
First, note that, since the morphism $T'\xrightarrow{f'} S$ has a section, it induces the splitting $T'=S\oplus T''=S\oplus \Ker (f')$ of generalized mixed Hodge structures. 
Now, we shall fix some notations. Let $\verb+{+c_{j}\verb+}+_{j=1}^{n}$ denote a basis of $T''$ over $\mathbb{Q}$ and put  
\begin{align}f_{|T''}\begin{pmatrix}
c_{1}\\
\vdots\\
c_{n}\\
\end{pmatrix}
=
\begin{pmatrix}
p_{1} & t_{1} & u_{1}\\
\vdots & \vdots & \vdots \\
p_{n} & t_{n} & u_{n}\\
\end{pmatrix}
\begin{pmatrix}
e\\
s_{1}\\
s_{2}\\
\end{pmatrix}
\end{align}
where $\verb+{+ p_{j}, t_{j}, u_{j}\verb+}+_{j=1}^{n}$ denotes elements of $\mathbb{Q}$. Then, it follows that $h(T'')$ is a $n$-dimensional subvector space of $V'$ spanned by the image $\verb+{+d_{j}=h(c_{j})\verb+}+_{j=1}^{n}$ over $\mathbb{Q}$. Thus, if we choose an element $v$ of $V'\verb+\+ h'(T'')$, the elements           $\verb+{+v,d_{j}\verb+}+_{j=1}^{n}$ form a basis of $V'$ over $\mathbb{Q}$. Define 
$$g(v)=p\alpha+q\beta \quad  \text{and} \quad g(d_{j})=p_{j}'\alpha+q_{j}'\beta.$$ 
Then, since we have $k\circ g=h'$ by the commutative diagram, we obtain $q\not=0$. Furthermore, since we also have $j\circ f=g\circ h$ by the commutative diagram and the action of $\verb+{+z_{D_{i}}\verb+}+_{i=1,2}$ on the image of $j$ in $V$ is given by $1$, it follows that we have $p_{j}'=p_{j}$ and $q_{j}'=0$ for all $1\leq j\leq n$. Thus, we can write   
\begin{align}g\begin{pmatrix}
v \\
d_{1}\\
\vdots\\
d_{n}\\
\end{pmatrix}
=
\begin{pmatrix}
p&q\\
p_{1}&0\\
\vdots&\vdots\\
p_{n}&0 
\end{pmatrix}
\begin{pmatrix}
\alpha\\
\beta\\
\end{pmatrix}
 \quad (q\not=0).
\end{align}
We shall show that we have $p_{j}=0$ ($1\leq j\leq n$). For simplicity, assume that the weight filtration $W_{j}T''$ (resp. $W_{j+1}T''$) of $T''$ is spanned by $\verb+{+c_{l+1},\cdots,c_{n}\verb+}+$ (resp. $\verb+{+c_{k},\cdots,c_{n}\verb+}+$) and that the quotient $\Gr_{j+1}^{W}T''$ has the Hodge structure of pure weight $0$. Then, by the argument of weights, it follows  that we have 
$$p_{j}=0 \quad(1\leq j\leq k-1 \  \text{and} \ l+1\leq j\leq n).$$
In particular, it follows from (3.2) that the image of $W_{j}T''$ under $f$ is contained in $i(S)$. By the commutative diagram, this means that we have $f_{|W_{j}T''}=0$. Thus, we obtain
$$t_{j}=0\quad \text{and}\quad u_{j}=0 \quad(l+1\leq j\leq n).$$
On the other hand, since $T''$ is an object of GMHS, there exist actions $\verb+{+x_{D_{i}}\verb+}+_{i=1,2}$ on $T''$    which are compatible with the actions of $\verb+{+w_{D_{i}}\verb+}+_{i=1,2}$ on $T$, that is, these satisfy $w_{D_{i}}f=fx_{D_{i}}$. Note that these actions induce involutions on $\Gr_{j+1}^{W}T''$ by definition.  
Let $\vec{c}$ be  the column vector $^{t}(c_{k},\cdots, c_{n})$ and put $x_{D_{i}}(\vec{c})=R_{i}\vec{c}$. Furthermore, $R_{i}'$ denotes the submatrix of $R_{i}$ which represents the residual action of $x_{D_{i}}$ modulo $W_{j}T''$.
Since we have $fx_{D_{i}}(\vec{c})=w_{D_{i}}f(\vec{c})$ and $p_{j}=t_{j}=u_{j}=0$ ($l+1\leq j\leq n$), it follows that we obtain
$$R_{i}'
\begin{pmatrix}
p_{k}e+t_{k}s_{1}+u_{k}s_{2}\\
\vdots\\
p_{l}e+t_{l}s_{1}+u_{l}s_{2}
\end{pmatrix}
=
\begin{pmatrix}
p_{k}e+(-p_{k}+t_{k})s_{1}+u_{k}s_{2}\\
\vdots\\
p_{l}e+(-p_{l}+t_{l})s_{1}+u_{l}s_{2}
\end{pmatrix}.
$$
If we put $\vec{p}=^{t}(p_{k},\cdots,p_{l})$ and $\vec{t}=^{t}(t_{k},\cdots,t_{l})$, this leads to $(R_{i}'-E)\vec{p}=0$ and $(R_{i}'-E)\vec{t}=-\vec{p}$. Since $R_{i}'$ denotes the matrix of the involution $x_{D_{i}}$ on $\Gr_{j+1}^{W}T''$, we have $(R_{i}'+E)(R_{i}'-E)=0$ and thus $(R_{i}'+E)\vec{p}=-(R_{i}'+E)(R_{i}'-E)\vec{t}=0$. Therefore, it follows that $\vec{p}$ is the zero-vector and that we obtain  
$$p_{j}=0 \quad (1\leq j\leq n).$$
Since $V'$ is also an object of GMHS, there exist actions $\verb+{+y_{D_{i}}\verb+}+_{i=1,2}$ on $V'$ which are compatible with the actions of $\verb+{+z_{D_{i}}\verb+}+_{i=1,2}$ on $V$, that is, these satisfy $z_{D_{i}}g=gy_{D_{i}}$. Put $y_{D_{i}}(v)=a_{0}v+\sum_{j=1}^{n}a_{j}d_{j}$. Since we have the formula (3.3) and $p_{j}=0$ for $1\leq j\leq n$, it follows that we obtain $gy_{D_{i}}(v)=a_{0}(p\alpha+q\beta)$. On the other hand,  since we have 
$z_{D_{i}}g(v)=z_{D_{i}}(p\alpha+q\beta)=p\alpha+q(c(D_{i})\alpha-\beta)$, the compatibility leads to 
$$a_{0}=-1\quad\text{and}\quad c(D_{i})=-\frac{2p}{q}\quad (i=1,2).$$
Here, note that we have $q\not=0$ by (3.3). This means that $c(D_{i})$ does not depend on $D_{i}$ and that this contradicts the assumption $c(D_{1})\not=c(D_{2})$. Thus, the Yoneda extension class $C(E)$ given by $E$ is non-trivial in $\Ext^{2}_{\text{GMHS}}(\mathbb{Q}_{M},S)$.  
\end{ex}

\end{document}